\documentclass[11pt]{amsart}
\usepackage{amssymb,amsmath}
\usepackage{amscd}
\usepackage{epsf}
\usepackage[all]{xy}
\usepackage{tikz}
\usepackage{mathrsfs}
\numberwithin{equation}{section}

\newtheorem{thm}{Theorem}[section]
\newtheorem{lem}[thm]{Lemma}

\newtheorem{defi}[thm]{Definition}

\newcommand{\inv}{^{-1}}

\newcommand{\C}{{\mathbb{C}}}
\newcommand{\Z}{{\mathbb{Z}}}
\newcommand{\R}{{\mathbb{R}}}
\newcommand{\N}{{\mathbb{N}}}

\newcommand{\Gr}{{\mbox{\rm Gr}}}

\newcommand{\Hom}{\mbox{\rm Hom}}

\newcommand{\trop}{\mbox{\rm trop}}
\newcommand{\Trop}{\mbox{\rm Trop}}

\newcommand{\Val}{\mbox{\rm Val}}

\newcommand{\cTG}{\mathcal{T}\!\Gr}

\newcommand{\an}{\mathrm{an}}
\newcommand{\degree}{\mathrm{deg}}
\newcommand{\slope}{\mathrm{slope}}
\begin{document}

\title{Analytification and tropicalization over non-archimedean fields}
\author{Annette Werner}
\date{}

\maketitle

\centerline{\bf Abstract:} In this paper, we provide an overview of recent progress on the interplay between tropical geometry and non-archimedean analytic geometry in the sense of Berkovich. After briefly discussing results by Baker, Payne and Rabinoff \cite{bpr} \cite{bpr2} in  the case of curves, we explain a result from \cite{chw} comparing  the tropical Grassmannian of planes to the analytic Grassmannian. We also give an overview of most of the results in \cite{grw}, where a general higher-dimensional theory is developed. In particular, we explain the construction of generalized skeleta in \cite{grw} which are polyhedral substructures of Berkovich spaces lending themselves to comparison with tropicalizations. We discuss the slope formula for the valuation of rational functions and explain two results on the comparison between polyhedral substructures of Berkovich spaces and tropicalizations. 
\small 
~\\[0.3cm]

\centerline{{\bf 2010 MSC: 14G22, 14T05} } 
%\centerline{{\bf Keywords: }}

\normalsize 
\section{Introduction}
Tropical varieties are polyhedral images of varieties over non-archime\-dean fields. They are obtained by applying the valuation map to a set of toric coordinates.  From the very beginning, analytic geometry was present in the systematic study of tropical varieties, e.g. in \cite{ekl} where rigid analytic varieties are used. The theory of Berkovich spaces which leads to spaces with nicer topological properties than rigid varieties is even better suited to study tropicalizations.
%For an overview and more explanation how rigid analytic geometry and Berkovich geometry are linked, see \cite{co}. 

A result of Payne \cite{pay} states that the Berkovich space associated to an algebraic variety is homeomorphic to the inverse limit of all tropicalizations in toric varieties. However, individual tropicalizations may fail to capture topological features of an analytic space. 
The present paper is an overview of recent results on the relationship between analytic spaces and tropicalizations. In particular, we address the question of whether a given tropicalization is contained in a Berkovich space as a combinatorial substructure.

A novel feature of Berkovich spaces compared to rigid analytic varieties is that they contain interesting piecewise linear combinatorial structures. In fact, Berkovich curves are, very roughly speaking, generalized graphs, where infinite ramifications along a dense set of points is allowed, see  \cite{ber}, chapter 4, \cite{baru} and \cite{bpr2}.

In higher dimensions, the structure of Berkovich analytic spaces is more involved, but still they often contain piecewise linear substructures as deformation retracts. This was a crucial tool in Berkovich's proof of local contractibility for smooth analytic spaces, see \cite{ber99} and \cite{ber04}. Berkovich constructs these piecewise linear substructures as  so-called skeleta of suitable models (or fibration of models). These skeleta basically capture the incidence structure of the irreducible components in the special fiber. 

In dimension one, Baker, Payne and Rabinoff studied the relationship between tropicalizations and subgraphs of Berkovich curves in \cite{bpr} and \cite{bpr2}. Their results show that every finite subgraph of a Berkovich curve admits a faithful, i.e. homeomorphic and isometric tropicalization. They also prove that every tropicalization with tropical multiplicity one everywhere is isometric to a subgraph of the Berkovich curve. 

As a first higher dimensional example, the Grassmannian of planes was studied in \cite{chw}. The tropical Grassmannian of planes has an interesting combinatorial structure and  is a moduli space for phylogenetic trees. It is shown in \cite{chw} that it is homeomorphic to a closed subset of the Berkovich analytic Grassmannian. 

In \cite{grw}, the higher dimensional situation is analyzed from a general point of view. This approach is based on a generalized notion of a Berkovich  skeleton which is associated to the datum of a semistable model plus a horizontal divisor. This  naturally leads to unbounded skeleta and generalizes a well-known construction on curves, see \cite{tyo} and \cite{bpr2}.

For every rational function $f$  with support in the fixed horizontal divisor it is shown in \cite{grw} that the ``tropicalization'' $\log |f|$ factors through a piecewise linear function on the generalized skeleton. This function satisfies a slope formula which is a kind of balancing condition around any 1-codimensional polyhedral face. Moreover, for every generalized skeleton there exists a faithful tropicalization -- where faithful in higher dimension refers to the preservation of the integral affine structures. In dimension one, this can be expressed via metrics. 
It is also proven that tropicalizations with tropical multiplicity one everywhere admit sections of the tropicalization map, which generalizes the above-mentioned result in \cite{bpr} on curves. 

The paper is organized as follows. In section 2 we collect basic facts on Berkovich spaces and tropicalizations, giving references to the literature for proofs and more details. In section 3 we briefly recall some of the results by Baker, Payne and Rabinoff \cite{bpr} on curves. Section 4 starts with the definition and basic properties of the tropical Grassmannian. Theorem \ref{thm:Grass} claims the existence of a continuous section to the tropicalization map on the projective tropical Grassmannian. We give a sketch of the proof for the dense torus orbit, where some constructions are easier to explain than in the general case. In section 5 we explain the construction of generalized skeleta from \cite{grw}, and in section 6 we investigate $\log|f|$ for rational functions $f$. In particular, Theorem \ref{slope_formula} states the slope formula. Section 7 explains the faithful tropicalization results in higher dimension.

{\bf Acknowledgements: }The author is very grateful to Walter Gubler for his helpful comments. 

\section{Berkovich spaces and tropicalizations}

\subsection{Notation and conventions} 
A non-archimedean field is a field with a non-archime\-dean absolute value. 
Our ground field is a non-archimedean field which is complete with respect to its absolute value.  Examples are the field $\mathbb{Q}_p$, which is the completion of $\mathbb{Q}$ after the $p$-adic absolute value, finite extensions of $\mathbb{Q}_p$ and also the $p$-adic cousin $\C_p$ of the complex numbers which is defined as the completion of the algebraic closure of $\mathbb{Q}_p$. The field of formal  Laurent series  $k((t))$ over an arbitrary base field $k$ is another example. Besides, we can endow any field $k$ with the trivial absolute value (which is one on all non-zero elements). A nice feature of Berkovich's general approach is that it also gives an interesting theory in the case of a trivially valued field. The reason behind this is that Berkovich geometry encompasses also points with values in transcendental field extensions -- and those may well carry interesting non-trivial valuations. 

Let $K$ be a complete non-archimedean field. We write $K^{\circ} = \{ x \in K: |x| \leq 1\}$ for the ring of integers in $K$, and $K^{\circ \circ} = \{ x´\in K: |x| < 1\}$ for the valuation ideal. The quotient $\tilde{K} = K^\circ / K^{\circ \circ}$ is the residue field of $K$. 
The valuation on $K^\times$ associated to the absolute value is given by $v(x) = - \log|x|$.
By $\Gamma = v(K^\times) \subset \R$ we denote the value group. 

A variety over $K$ is an irreducible, reduced and separated scheme of finite type over $K$.

\subsection{Berkovich spaces}\label{section:berkovich} Let us briefly recall some basic results about Berkovich spaces. This theory was developed in the ground-breaking treatise \cite{ber}. The survey papers
\cite{co}, \cite{du} and \cite{tem} provide additional information. For background information on non-archimedean fields and Banach algebras and for an account of rigid analytic geometry see \cite{bgr} and  \cite{bo}. 

For  $n \in \mathbb{N}$ und every $n-$tuple $r = (r_1, \ldots, r_n)$ of positive real numbers we define the associated \emph{Tate algebra} as 
\[K\{ r\inv x\} =  \{\sum_{I = (i_1, \ldots, i_n) \in \N_0^n} a_I x^I: |a_I| r^I \rightarrow 0 \mbox{ as }|I| \rightarrow \infty\}.\]
Here we put $x=(x_1, \ldots, x_n)$, and for $I = (i_1, \ldots, i_n)$ we write $x^I = x^{i_1} \ldots x^{i_n}$. Moreover, we define $|I| = i_1 + \ldots +  i_n$. 

A Banach algebra $A$  is called \emph{$K$-affinoid}, if there exists a surjective $K$-algebra homomorphism $\alpha: K\{r\inv x\} \rightarrow A$ for some $n$ and $r$ such that the Banach norm on $A$ is equivalent the quotient seminorm induced by $\alpha$. 
If such an epimorphism can be found with  $r = (1, \ldots, 1)$, then $A$ is called \emph{strictly $K$-affinoid}. In rigid analytic geometry only strictly $K$-affinoid algebras are considered (and called affinoid algebras).

The \emph{Berkovich spectrum} $\mathcal{M}(A)$ of an affinoid algebras is defined as the set of all  bounded (by the Banach norm), multiplicative seminorms on $A$. It is endowed with the coarsest topology such that all evaluation maps on functions in $A$ are continuous. 

If $x$ is a seminorm on an algebra $A$, and $f \in A$, we follow the usual notational convention and write $|f(x)|$ for $x(f)$, i.e. for the real number which we get by evaluating the seminorm $x$ on $f$. 

The \emph{Shilov boundary} of a Berkovich spectrum $\mathcal{M}(A)$ of a $K$-affinoid algebra $A$ is the unique minimal subset $\Gamma$ (with respect to inclusion) such that for every $f \in A$ the evaluation map $\mathcal{M}(A) \rightarrow \R_{\geq 0}$ given by $x \mapsto |f(x)|$ attains its maximum on $\Gamma$. Hence for every $f \in A$ there exists a point $z$ in the Shilov boundary such that $|f(z)| \geq |f(x)|$ for all $x \in \mathcal{M}(A)$. The Shilov boundary of $\mathcal{M}(A)$ exists and is a finite set by \cite{ber}, Corollary 2.4.5.

The Berkovich spectrum of a $K$-affinoid algebra carries a sheaf of analytic functions which we will not define here. Some care is needed since, very roughly speaking, not all coverings are suitable for glueing analytic functions.
Analytic spaces over $K$ are ringed spaces which are locally modelled on Berkovich spectra of affinoid algebras. 

There is a \emph{GAGA-functor}, associating to every $K$-scheme $X$ locally of finite type an analytic space $X^{\an}$ over $K$. It has the property that $X$ is connected, separated over $K$ or proper over $K$, repectively, if and only if the topological space $X^{\an}$ is arcwise connected, Hausdorff, or compact, respectively. If $X = \mbox{Spec} R$ is affine, the topological space $X^{\an}$ can be identifed with the set of all multiplicative seminorms on the coordinate ring $R$ extending the absolute value on $K$. This space is endowed with the coarsest topology such that for all $f \in R$ the evaluation map on $f$ is continuous. 

With the help of the GAGA-functor we associate analytic spaces to algebraic varieties. 
Another way of obtaining analytic spaces is via \emph{admissible formal schemes}. An admissible formal scheme over $K^\circ$ is, roughly speaking, a formal scheme over $K^\circ$  such that the formal affine building blocks are given by $K^\circ$-flat algebras of the form $K^\circ\{x_1, \ldots, x_n\} / \mathfrak{a}$ for a finitely generated ideal $\mathfrak{a}$ in the ring $K^\circ\{x_1, \ldots, x_n\} = \{\sum_{I \in \mathbb{N}_0^n} a_I x^I: a_I \in K^\circ \mbox{ and }|a_I| \rightarrow 0 \mbox{ as }|I| \rightarrow \infty\}$. For a precise definition see \cite{bo} or \cite{co}.

An admissible formal scheme $\mathcal{X}$ has a natural \emph{analytic generic fiber} $\mathcal{X}_\eta$. On the formal affine building block given by $K^\circ\{x_1, \ldots, x_n\} / \mathfrak{a}$ the analytic generic fiber is the Berkovich spectrum of $K \{x_1, \ldots, x_n\} / \mathfrak{a}K \{x_1, \ldots, x_n\} $.  

If we start with a $K^\circ$- scheme $\mathscr{X}$, which is separated, flat and of finite presentation, its completion with respect to any element $\pi \in K^{\circ  \circ} \backslash \{0\}$ is an admissible formal scheme $\mathcal{X}$. This  has an analytic generic fiber $\mathcal{X}_\eta$. On the other hand, we can apply the GAGA-functor to the algebraic generic fiber $X = \mathscr{X} \otimes_{K^\circ} K$ and get another analytic space $X^\an$. They are connected by a morphism $\mathcal{X}_\eta \hookrightarrow X^\an$, which is an isomorphism if $\mathscr{X}$ is proper over $K^\circ$. In the basic example $\mathscr{X} = \mbox{Spec}K^\circ [x]$ we get the natural inclusion $\mathcal{M}(K\{x\}) \hookrightarrow (\mathbb{A}^1_K)^\an$, whose image is the set of all multiplicative seminorms on the polynomial ring $K[x]$ which extend the absolute value on $K$ and whose value on $x$ is bounded by $1$. 

\subsection{Tropicalization}\label{section:tropicalization}

We start by considering a split torus $T = \mathrm{Spec}K[x_1^{\pm 1}, \ldots, x_n^{\pm 1}]$. 
Evaluating seminorms on characters we get a map 
$\trop: T^\an \rightarrow \R^n$ on the associated Berkovich space, which is given by
\[\trop(p) = (-\log|x_1(p)|, \ldots, -\log|x_n(p)|).\]
Here we are using the valuation map to define tropicalizations. In some papers, tropicalizations are defined via the negative valuation map, i.e. the logarithmic absolute value of the coordinates. 

If $X$ is a variety over $K$ together with a closed embedding $\varphi: X \hookrightarrow T$, we consider the composition 
\begin{equation}
\label{eq:trop}
\trop_{\varphi}: X^\an \stackrel{\varphi^\an}{\longrightarrow } T^\an \stackrel{\tiny \trop}{\longrightarrow} \R^n
\end{equation}
 of the tropicalization map with the embedding $\varphi$. 
 Note that the map $\trop_\varphi$ is continuous.
Its image $\Trop_\varphi(X) = \trop_\varphi (X^\an)$ is the support of a polyhedral complex $\Sigma$ in $\R^n$. This complex is integral $\Gamma$-affine in the sense of section \ref{section:integral}. It has pure dimension $d = \dim(X)$. 
For a nice introductory text on tropical geometry see the textbook \cite{ms}. The survey paper \cite{gu13} is also very useful.

Note that for all $\omega \in \Trop_\varphi(X)$ the preimage $\trop_\varphi^{-1}(\{\omega\})$ can be identified with the Berkovich spectrum of an affinoid algebra, see \cite{gu07}, Proposition 4.1.

For every point $\omega \in \Trop_\varphi(X)$ there is an associated \emph{initial degeneration}, which is the special fiber of a model of $X$ over the valuation ring in a suitable non-archimedean extension field of $K$. The geometric number of irreducible components of this special fiber is the \emph{tropical multiplicity} of $\omega$. There is a balancing formula for tropicalizations involving tropical multiplicities on all maximal-dimensional faces of $\Trop_\varphi(X)$ around a fixed codimension one face, see \cite{ms}, chapter 3.4.

Let $\overline{K}$ be the algebraic closure of $K$. It can be endowed with an absolute value extending the one on $K$. By means of the coordinates $x_1, \ldots, x_n$ we can identify $T(\overline{K}) = \overline{K}^n$. Then we can define a natural tropicalization map $\trop: T(\overline{K}) \rightarrow \R^n$ by $\trop(t_1, \ldots, t_n) = (-\log|t_1|,\ldots, -\log|t_n|)$. Note that every point $t = (t_1, \ldots, t_n)$ in $T(\overline{K})$ gives rise to the multiplicative seminorm $f \mapsto |f(t)|$ on the coordinate ring $K[x_1^{\pm 1}, \ldots, x_n^{\pm 1}]$, which is an element in $T^\an$. Hence the tropicalization map  on $T(\overline{K})$ is induced by the map $\trop$ on $T^\an$.

We could define the tropicalization $\Trop_\varphi(X)$ without recourse to Berkovich spaces. If the absolute value on $K$ is non-trivial, the tropicalization $\Trop_\varphi(X)$ is  equal to the closure of the image of the map
\[X(\overline{K}) \stackrel{\varphi}{\longrightarrow} T(\overline{K}) \stackrel{\tiny \trop}{\longrightarrow} \R^n.\]
Considering a tropical variety as the image of an analytic space makes some topological considerations easier. For example, in this way it is evident that the tropicalization of a connected variety is connected since it is a continuous image of a connected space. 

Let $Y$ be a toric variety associated to the fan $\Delta$ in $N_\mathbb{R}$, where $N$ is the cocharacter group of the dense torus $T$. Then $\Delta$ defines a natural partial compactification $N_\R^\Delta$ (due to Kashiwara and Payne)   of the space $N_\R$. Roughly speaking, we compactify cones by dual spaces. For details see \cite{pay}, section 3.
There is also a natural  tropicalization map
\[\trop: Y^{\an} \rightarrow N_\R^\Delta,\]
which extends the tropicalization map on the dense torus.
As an example, we consider $Y = \mathbb{P}^n_K$ with its dense torus $\mathbb{G}_{m,K}^n$ and $N = \Z^n$. Put $\overline{\R} = \R \cup \{\infty\}$ and endow it with the natural topology such that half-open intervalls $]a, \infty]$ form a basis of open neighbourhoods of $\infty$.  Then the associated compactification of $N_\R = \R^n$ is the tropical projective space
\[\mathbb{TP}^n = \big(\overline{\R}^{n+1}
\smallsetminus \{(\infty, \ldots, \infty)\}\big) /
\R (1, \ldots, 1)
  \]
which is endowed with the product-quotient topology. 
The Berkovich projective space $(\mathbb{P}^n_K)^\an$ can be described as the set of equivalence classes of elements in $(\mathbb{A}^{n+1}_K)^\an \backslash \{0\}$ with respect to the following equivalence relation: Let $x$ and $y$ be points in $(\mathbb{A}^{n+1}_K)^\an \backslash \{0\}$, i.e.  multiplicative seminorms on $K[x_0, \ldots, x_n]$ extending the absolute value on $K$. They are equivalent if and only if there exists a constant $c >0$ such that for every homogeneous polynomial $f$ of degree $d$ we have $|f(y)| = c^d |f(x)|$. 

We can describe the tropicalization map on the toric variety $\mathbb{P}^n_K$ as follows:
\begin{equation}\label{eq:tropical_projective}
\trop: (\mathbb{P}^n_K)^\an \rightarrow \mathbb{TP}^n
\end{equation}
maps the class of the seminorm $p$ on $K[x_0, \ldots, x_n]$ to the class
$(-\log|x_0(p)|, \ldots, -\log|x_n(p)|) + \R(1, \ldots, 1)$ in tropical projective space. 

An important result about the relation between tropical and analytic geometry due to Payne says that for every quasi-projective variety $X$ over $K$ the associated analytic space $X^{\an}$ is homeomorphic to the inverse limit over all $\Trop_\varphi(X)$, where the limit runs over all all closed embeddings $\varphi: X \hookrightarrow Y$ in a quasi-projective toric variety $Y$ (see \cite{pay}, Theorem 4.2). For a generalization which omits the quasi-projectivity hypothesis see \cite{fgp}.

Hence tropicalizations are combinatorial images of analytic spaces, which recover the full analytic space in the projective limit. Individual tropicalizations however may not faithfully depict all topological features of the analytic space.
We will discuss a basic example at the beginning of the next section.

%Using Payne's extended tropicalizations of toric varieties \cite{pay}, we can generalize this construction to obtain ``standard skeleta´´ of toric varieties. Note that they depend on the choice of torus coordinates, i.e. on the choice of a basis of the character group.  

\section{The case of curves} \label{section:curves}
Let $K$ be a field which is algebraically closed and complete with respect to a non-archimedean, non-trivial absolute value, and let $X$ be a smooth curve over $K$.

As an example, let us consider an elliptic curve $E$ in Weierstrass form 
$y^2 = x^3 + ax +b$ over $K$. The affine curve $E$ has a natural embedding in $\mathbb{A}^2_K$. The intersection $E_0$ of $E$ with $\mathbb{G}_{m,K}^2$ has a natural tropicalization in $\R^2$, which is one of the trees depicted in the next figure.

\begin{figure}[htf]

  \begin{minipage}[l]{0.45\linewidth}
    \includegraphics[height=3cm]{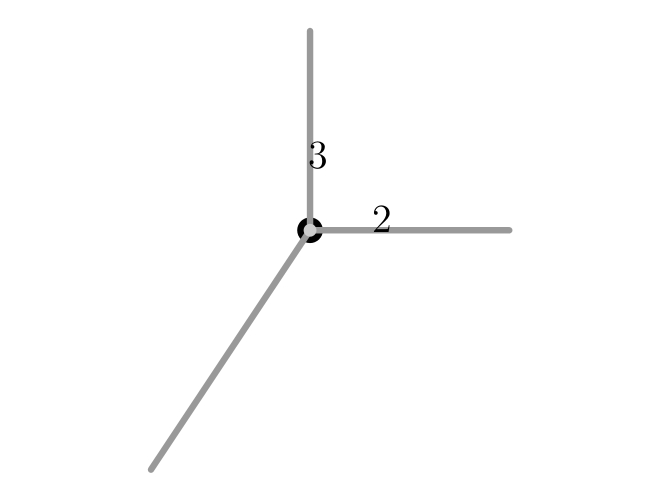}
    \end{minipage}
    \begin{minipage}[r]{0.3\linewidth}
    \includegraphics[height=4cm]{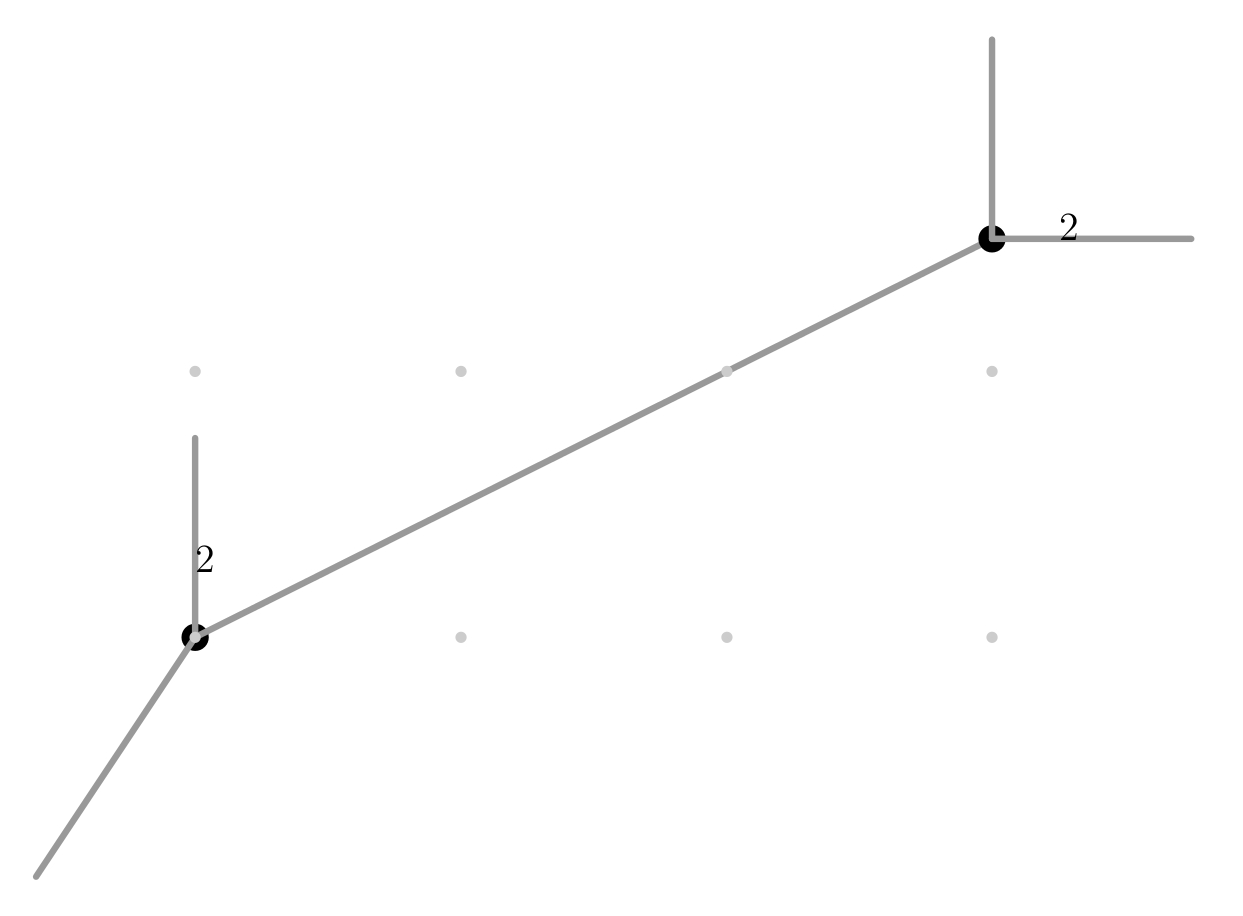}
    \end{minipage}
    
 \caption{Tropicalization of an elliptic curve in Weierstra{\ss} form. The left hand side shows the case $3v(a) \geq 2 v(b) $, the right hand side the case $3 v(a) < 2 v(b)$. The numbers on the edges indicate the tropical multiplicities. If no number is given, the multiplicity is one.}
 \label{figure:weier}
    
    \end{figure}

    If however $E$ is a Tate curve, then the analytic space $E^{an}$ contains a circle, i.e. it has topological genus $1$. Hence the Weierstrass tropicalization does not faithfully depict the topology of the analytic space.

As stated before, if $X$ is a smooth curve over $K$, the Berkovich space $X^{an}$ is some kind of generalized graph, allowing infinite ramification along a dense set of points. In particular, it makes sense to talk about its  leaves. The   space of non-leaves $H_0(X^\an)$ in the Berkovich analytification $X^\an$ admits a natural metric which is defined via semistable models, see \cite{bpr2}, section 5.3. 
Any tropicalization of $X$ with respect to a closed embedding in a torus carries a natural metric which is locally given by lattice length on each edge (and globally by shortest paths). 
 Another problem in the case of the Weierstra{\ss} tropicalizations from figure \ref{figure:weier} is the fact that tropical multiplicites $\neq 1$ are present. This may indicate that the tropicalization map will not be isometric on a subgraph, see  \cite{bpr}, Corollary 5.9.
    
It is explained in  \cite{bpr}, Theorem 6.2, how one can construct better tropicalizations of the Tate curve. 
    
Based on a detailed study of the structure of analytic curves 
%and a crucial result comparing the affinoid algebra given by the preimage of tropicalization to the affinoid algebra associated to the analytic generic fiber of the initial degeneration (\cite{bpr}, section 4)
there are the following two important comparison theorems in \cite{bpr}.

\begin{thm}[\cite{bpr}, Theorem 5.20] Let $\Gamma$ be any finite subgraph of $X^\an$ for a smooth, complete curve $X$ over $K$. Then there exists a closed immersion $\varphi: X \hookrightarrow Y$ into a toric variety $Y$ with dense torus $T$ such that the restriction $\varphi_0: X_0 = X \cap \varphi^{-1}(T) \hookrightarrow T$ induces a tropicalization map $\trop_{\varphi_0}: X_0^\an \rightarrow \R^n$ which maps $\Gamma$ homeomorphically and isometrically onto its image. 
\end{thm}

\begin{thm}[\cite{bpr}, Theorem 5.24] Consider a smooth curve $X_0$ over $K$ and a closed immersion $\varphi: X_0 \hookrightarrow T$ into a split torus with associated tropical variety $\Trop_{\varphi_0}(X_0)$. If $\Gamma'$ is a compact connected subset of $\Trop_{\varphi_0}(X_0)$ which has tropical multiplicity one everywhere, then there exists a unique closed subset $\Gamma$ in $H_0(X_0^\an)$ mapping homeomorphically onto $\Gamma'$, and this homeomorphism is in fact an isometry.
\end{thm}

We will discuss higher-dimensional generalizations of these theorems in section \ref{section:faithful}.

\section{Tropical Grassmannians} \label{section:Grass}

\subsection{The setting}\label{tropgrass} In view of the results  \cite{bpr} for curves, it is a natural question whether they can be generalized to varieties of higher dimensions. As a first example, the tropical Grassmannian of planes was studied in \cite{chw}. In this section we discuss the main theorem of this paper.

For natural numbers $d \leq n$ we denote by $\Gr(d,n)$ the Grassmannian of $d$-dimensional subspaces of $n$-space. 
The tropical Grassmannian is defined as the tropicalization of $\Gr(d,n)$ with respect to the Pl\"ucker embedding $\varphi: \Gr(d,n) \rightarrow \mathbb{P}^{\binom{n}{d} -1}$. Recall that the Pl\"ucker embedding maps the point corresponding to the $d$-dimensional subspace $W$ to the line in $\binom{n}{d}$-space given by the $d$-th exterior power of $W$. 

During this section we will deviate from the exposition in section \ref{section:tropicalization} and  consider tropicalizations with respect to the negative valuation map. The choice between $(\mbox{min},+)$ and $(\mbox{max},+)$ tropical geometry is always a difficult one. 

Hence our tropical projective space is
\[\mathbb{T}\mathbb{P}^{\binom{n}{d} -1} = \big( (\R \cup \{- \infty\})^{\binom{n}{d}}
\smallsetminus \{(-\infty, \ldots, -\infty)\}\big) /
\R (1, \ldots, 1), 
  \]
and the tropical Grassmannian $\mathcal{T}\Gr(d,n) = \mbox{Trop}_\varphi(\Gr(d,n))$ is defined as the image of 
$\trop \circ \varphi^\an: \Gr(d,n)^\an \rightarrow \mathbb{T}\mathbb{P}^{\binom{n}{d} -1}$, where  $\trop$ is given by the map
\[ p \mapsto (\log|x_0(p)|, \ldots, \log|x_{\binom{n}{d}}(p)|) + \R(1, \ldots, 1)\]
on the analytic projective space, see
 (\ref{eq:tropical_projective}).

Tropical Grassmannians were first studied by Speyer and Sturmfels \cite{ss} who focused on the toric part $\Gr_0(d,n) = \Gr(d,n) \cap \varphi\inv \big(\mathbb{G}_{m,K}^{\binom{n}{d}-1}\big)$ which is embedded in the torus  $\mathbb{G}_{m,K}^{\binom{n}{d}-1}$ via $\varphi$. The tropicalization $\mathcal{T}\Gr_0(d,n) = \Trop_\varphi(\Gr_0(d,n))$ (which is called $\mathcal{G}'_{d,n}$ in \cite{ss}) is a fan of dimension $d (n-d)$ containing a linear space of dimension $n-1$, see \cite{ss}, section 3.

Moreover, if $d=2$, Speyer and Sturmfels proved that  $\mathcal{T} \Gr_0(2,n)$ can be identified with the space of phylogenetic trees, see \cite{ss}, section 4. For our purposes, a phylogenetic tree on $n$ leaves
is a pair $(T,\omega)$, where  $T$ is a finite combinatorial tree with no degree-two vertices together with a
labeling of its leaves in bijection with $\{1, \ldots,n \}$, and $\omega$ is a real-valued function on the set of edges of $T$.  The tree $T$ is called the combinatorial type of
the phylogenetic tree $(T, \omega)$. For every phylogenetic tree $(T, \omega)$ and all $i \neq j$ in $\{1, \ldots, n\}$ we denote by $x_{ij}$ the sum of the weights along the uniquely determined path from leaf $i$ to leaf $j$. 
This tree-distance function satisfies the four-point-condition, which states that for all pairwise distinct indices $i,j,k,l$ in $\{1, \ldots, n\}$ the maximum among 
\[x_{ij} + x_{kl}, \quad x_{ik} + x_{jl},\quad x_{il} + x_{jk}\]
is attained at least twice. 

\subsection{A section of the tropicalization map}\label{sectionforgrass}
For the rest of this section we will always consider the case $d=2$. 
In \cite{chw}, we investigate the full projective Grassmannian $\mathcal{T} \Gr(2,n) = \Trop_\varphi \Gr(2,n)$ in tropical projective space $\mathbb{T}\mathbb{P}^{\binom{n}{2} -1}$. 
The main result of this paper is the following theorem: \begin{thm}[\cite{chw}, Theorem 1.1]\label{thm:Grass}
  There exists a continuous section $\sigma\colon \cTG(2,n)\to
  \Gr(2,n)^{\an}$ of the tropicalization map $\trop \circ \varphi^\an \colon
  \Gr(2,n)^{\an} \to \cTG(2,n)$. Hence, the tropical Grassmannian
  $\cTG(2,n)$ is homeomorphic to a closed subset of the Berkovich
  analytic space $\Gr(2,n)^{\an}$.
\end{thm}

Note that we are not considering semistable models or their Berkovich skeleta in this approach. 

The first idea one might have is to look at the big open cells of the Grassmannians. Since $d =2$, the coordinates of the ambient projective space
are indexed by the two element subsets $\{i,j\}$ of $\{1, \ldots, n\}$. If $i<j$, we write $p_{ij}$ for this coordinate, and we put $p_{ji} = - p_{ij}$. We intersect the Grassmannian with the standard open affine covering of projective space and get open affine subvarieties 
\[ U_{ij} = \Gr(2,n) \cap \varphi\inv \{p_{ij} \neq 0\}. \]
The affine Pl\"ucker coordinates are $u_{kl} = p_{kl} / p_{ij}$. 
Since the Pl\"ucker ideal is generated by the relations $p_{ij} p_{kl} - p_{ik} p_{jl} + p_{il} p_{jk} = 0$ for $\{i,j,k,l\}$ running over the four-element subsets of $\{1, \ldots, n\}$,  the subvariety $U_{ij}$ can be identified with $\mathbb{A}_K^{n(n-2)}$ by means of  the coordinates $u_{ik} = p_{ik} / p_{ij}$ and $u_{jk} = p_{jk} / p_{ij}$ for all $k$ different from $i$ and $j$. 

Recall that $(\mathbb{A}_K^{n(n-2)})^\an$ is the set of all multiplicative seminorms on $K[u_{ik}, u_{jk}: k \notin \{i,j\}]$ which extend the absolute value on $K$. As explained in section 2, the toric variety $\mathbb{A}_K^{n(n-2)}$ has a natural tropicalization. With the sign conventions in the present section, it is given by
\begin{eqnarray*}
\trop: (\mathbb{A}_K^{n(n-2)})^\an & \rightarrow &  (\R \cup \{- \infty\})^{2(n-2)}\\
p & \mapsto & (\log|u_{ik}(p)|, \log|u_{jk}(p)|)_{k \notin \{i,j\}}
\end{eqnarray*}
This induces the tropicalization map 
\[
\trop: U_{ij}^\an \longrightarrow (\mathbb{R}\cup \{- \infty\})^{2(n-2)}.\]

Now the tropicalization map on an affine space $(\mathbb{A}^N_K)^\an$ has a natural section, which is given by the map 
$\delta: (\mathbb{R} \cup \{- \infty\})^N \rightarrow (\mathbb{A}^N_K)^\an$, mapping a point $r=(r_1, \ldots, r_N) \in (\R \cup \{- \infty\})^N$ to the seminorm 
\begin{equation}
\label{eq:skeleton}
\delta(r)\colon K[x_1,\ldots, x_N] \longrightarrow \R_{\geq 0} \qquad
  |\sum_{\alpha \in \mathbb{N}_0^N} c_{\alpha}{x}^{\alpha} (\delta(r))|
 =\max_\alpha \{ |c_\alpha| \prod_{i=1}^N
 \exp({r_i \alpha_i})\}.
 \end{equation}
Here we put $\exp((-\infty)\alpha) = 0$ for all $\alpha \in \mathbb{N}_0$. Note that the restriction of $\delta$ to $\mathbb{R}^N$ is a map from $\mathbb{R}^N$ to the torus $(\mathbb{G}_m^{N})^\an$. We call the image of $\delta$ on $(\R \cup \{-\infty\})^N$ the standard skeleton of $(\mathbb{A}_K^{N})^\an$, and the image of $\delta|_{\R^N}$ the standard skeleton of the torus $ (\mathbb{G}_m^{N})^\an$.

However, the map $\delta$ for $N = n(n-2)$ does not in general provide a section of the tropicalization map on the whole of $U_{ij}$.  Consider $n\geq 4$, and let $x = (x_{kl})_{kl}$ be a point in the tropical Grassmannian $\cTG(2,n)$ which lies in the image of $U_{ij}$ under the tropicalization map.  Let $\omega$ be the projection to the affine coordinates $\omega= (x_{ik}-x_{ij}, x_{jk}-x_{ij})_{k \notin\{i,j\}} \in \R^{2(n-2)}$. Then $\delta(\omega)$ is a point in the Berkovich space $U_{ij}^\an$ with $\log |u_{ik}(\delta(\omega))| = x_{ik}-x_{ij}$ and $\log|u_{jk}(\delta(\omega))| = x_{jk}-x_{ij}$. Since $u_{kl} := p_{kl} / p_{ij} = u_{ik} u_{jl} -  u_{il} u_{jk}$ by the Pl\"ucker relations, the definition of $\delta$ gives us $ \log|u_{kl}(\delta(\omega))| = \max\{ x_{ik} + x_{jl}, x_{il} + x_{jk}\}- 2 x_{ij}$. Hence $\delta$ can only provide a section of the tropicalization map if this maximum is equal to $x_{kl}+x_{ij}$. This may fail if the labelled tree $T$ has the wrong shape, e.g. if it looks like this:
\begin{figure}[h!]
  \centering
  \begin{minipage}[c]{0.45\linewidth}
    \includegraphics[scale=0.45]{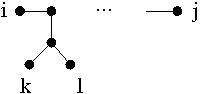}
  \end{minipage}
\end{figure}

The strategy for the proof of Theorem \ref{thm:Grass} is to compare the tropical Grassmannian with a standard skeleton of an affine space on a smaller piece of the tropicalization, namely on the part consisting of phylogenetic trees such that the underlying combinatorial tree has the right shape. If we want to make  this precise, the definition of these smaller pieces is quite involved. A considerable part of the difficulties is due to the fact that we  take into account the boundary strata of the Grassmannian in projective space. In order to explain the general strategy we will from now on restrict our attention to the torus part of the tropical Grassmannian. The general case can be found in \cite{chw}.

\subsection{Sketch of proof in the dense torus orbit}\label{sectiontoruspart} 
In this section we  explain the proof of Theorem \ref{thm:Grass} for the subset $\cTG_0(2,n)$ of the tropical Grassmannian $\cTG(2,n)$. Recall that $\cTG_0(2,n)$ is the tropicalization of the dense open subset $\Gr_0(2,n)$ of the Grassmannian which is mapped to the torus via the Pl\"ucker map. 

We fix a pair $ij$ as above and work in the big open cell $U_{ij}$. The coordinate ring $R_{ij}$ of $U_{ij}$ is a  polynomial ring  $K[u_{ik}, u_{jk}: k \notin \{i,j\}]$ in $2(n-2)$ variables. 
The other Pl\"ucker coordinates are expressed as $u_{kl} = u_{ik} u_{jl} - u_{il} u_{jk}$ in $R_{ij}$. 
Note that the affine variety $\Gr_0(2,n)$ is contained in $U_{ij}$. The coordinate ring of $\Gr_0(2,n)$ is equal to the localization of $R_{ij}$ after the multiplicative subset generated by all $u_{kl}$ for $\{k,l\} \neq \{i,j\}$. 

We also fix a labelled tree $T$ with $n$ leaves $1, \ldots, n$, and arrange $T$ as in the following figure

 \begin{figure}[h!]
  \centering
  \begin{minipage}[c]{0.45\linewidth}
    \includegraphics[scale=0.45]{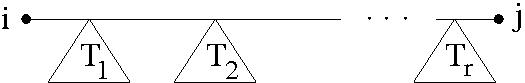}
  \end{minipage}
  \caption{Arrangement of $T$}
  \label{tree}
\end{figure}

with subtrees $T_1, \ldots, T_r$. 

\begin{defi}\label{def:orderWithCherryProperty} Let  $\preceq$ be a partial order on the set
  $\{1, \ldots, n\} \smallsetminus \{i,j\}$. We write $k \prec l$ if $k \preceq l$ and $k \neq l$. Then $\preceq$ has the \emph{cherry property} on $T$ with respect to
  $i$ and $j$ if the following conditions hold:
  
\begin{enumerate}
\item[(i)] Two leaves of different subtrees $T_a$ and $T_b$ for $a,b \in \{1, \ldots r\}$ as in figure \ref{tree} cannot be
  compared by $\preceq$.
\item[(ii)] The partial order $\preceq$ restricts to a total order on the leaf
  set of each $T_a$, $a=1,\ldots, r$.
\item[(iii)] If $k\prec l\prec m$, then either $\{k,l\}$ or $\{l,m\}$ is a
  \emph{cherry of the quartet} $\{i,k,l,m\}$, i.e. the subtree given by this quartet of leaves contains a node which is adjacent to both elements of the cherry:
  \begin{figure}[hbf]
  \centering
  \begin{minipage}[c]{0.45\linewidth}
    \includegraphics[scale=0.45]{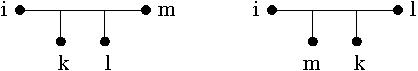}
  \end{minipage}
  \caption{Cherry property}
\end{figure}
\end{enumerate}
  \end{defi}

 An induction argument  shows the following lemma:
\begin{lem}[\cite{chw}, Lemma 4.7]\label{lm:order} Fix a pair of indices $i,j$, and let $T$
  be a tree on $n$ labelled leaves. Then, there exists a partial order
  $\preceq$ on the set $\{i, \ldots, n\} \smallsetminus \{i,j\}$ that has the cherry
  property on $T$ with respect to $i$ and $j$.
\end{lem}
For an example of the inductive construction of such a partial order see \cite{chw}, Figure 3. 

The leaves in each subtree $T_a$  are totally ordered, say as $s_1 \prec s_2 \prec \ldots \prec s_p$, if $T_a$ contains $p = p(a)$ leaves. We consider the variable set 
$I_a = \{ u_{i s_1}, \ldots, u_{i s_p}\} \cup \{u_{j s_1}, u_{s_1 s_2}, \ldots u_{s_{p-1} s_p}\}$. 
Then $I = I_1 \cup \ldots \cup I_r$ is a set of $2(n-2)$ affine Pl\"ucker coordinates of the form $u_{kl} \in R_{ij}$. 

We can successively reconstruct the variables of the form $u_{jl}$ which are not contained in  $I$ as follows: In the tree $T_a$ with leaves $s_1 \prec s_2 \prec \ldots \prec s_p$, we have $u_{s_1 s_2} = u_{i s_1} u_{j s_2} - u_{i s_2} u_{j s_1}$, hence
\[u_{j s_2} = u_{i s_1}\inv ( u_{s_1 s_2} + u_{i s_2} u_{j s_1}).\]
The right hand side  is an expression in the variables contained in the coordinate set $I$ with $u_{i s_1}$ inverted. 
Now we use the relation $u_{s_2 s_3} = u_{i s_2} u_{j s_3} - u_{i s_3} u_{j s_2}$ to express $u_{j s_3} = u_{i s_2}\inv (u_{s_2 s_3} + u_{i s_3} u_{j s_2})$. Plugging in the expression for $u_{j s_2}$ we can write $u_{j s_3}$ as a polynomial in the variables in $I$ plus all $u_{ik}\inv$. 

Proceeding by induction,  we find for all $m \neq i,j$ that 
\[ u_{jm} \in K[u_{kl}: u_{kl} \in I][u_{ik}\inv : k \neq i,j]\]
and hence
\begin{equation}\label{eq:birational}
K[u_{kl}: u_{kl }\in I] \subset R_{ij} \subset K[u_{kl}: u_{kl} \in I][u_{ik}\inv : k \neq i,j].
\end{equation}
This shows that the variable set $I$  generates the function field $\mbox{Quot}(R_{ij})$ of $U_{ij}$. 
Recall that the coordinate ring of $\Gr_0(2,n)$ is  $K[\Gr_0(2,n)] = S^{-1} R_{ij}$,
where $S$ is the multiplicative subset of $R_{ij}$ generated by all $u_{kl}$ for $\{kl\} \neq \{ij\}$. By the previous result, $K[\Gr_0(2,n)]$  is equal to the localization of the Laurent polynomial ring $K[u_{kl}^{\pm}: u_{kl} \in I]$ after the multiplicative subset generated by \emph{all} $u_{kl}$ expressed as Laurent polynomials in the coordinates contained in $I$.  

\begin{defi} Let $\mathcal{C}_T$ be the cone in  $\cTG_0(2,n)$ whose interior corresponds to the phylogenetic trees with underlying tree $T$.
Let $x= (x_{kl})_{kl} + \mathbb{R}(1, \ldots, 1)$ be a point in $\mathcal{C}_T \subset \cTG_0(2,n)$. We associate to it a point $\sigma_T^{ij}(x)$ in $\Gr_0(2,n)^\an$, i.e. a multiplicative seminorm on the coordinate ring of $\Gr_0(2,n)$, as follows. For every Laurent polynomial $f = \sum_{\alpha } c_\alpha u^\alpha \in K[u_{kl}^\pm : u_{kl} \in I]$ (where $\alpha$ runs over $\mathbb{Z}^I$) we put
\[|f(\sigma_T^{ij}(x))| = \max_\alpha \big\{|c_\alpha| \prod_{u_{kl} \in I} \exp\big( \alpha_{kl} (x_{kl} - x_{ij})\big) \big\}, \quad \mbox{ where } \alpha = (\alpha_{kl})_{u_{kl} \in I}.\]
This defines a multiplicative norm on $K[u_{kl}^\pm: \{kl\} \in I]$ which has a  unique extension  to a multiplicative norm on the localization $K[\Gr_0(2,n)]$. Let $\sigma_T^{ij}(x)$ be the resulting point in the Berkovich space $\Gr_0(2,n)^\an$.
\end{defi}

Since for every $f \in K[\Gr_0(2,n)]$ the evaluation map on $\mathcal{C}_T$ given by
\[x \mapsto |f(\sigma_T^{ij}(x))|\]
is continuous, we have constructed a continuous map
\[\sigma_T^{ij}: \mathcal{C}_T \rightarrow \Gr_0(2,n)^\an.\]

We want to show that it is a section of the tropicalization map 
 $\trop \circ \varphi^\an : \, \Gr_0(2,n)^\an \rightarrow \cTG_0(2,n)$ on $\trop\inv(\mathcal{C}_T)$, which amounts to checking that 
 \[\log |u_{kl}(\sigma_T^{ij}(x))| =  x_{kl} - x_{ij} \mbox{ for all } \{kl\}\neq\{ij\} \mbox{ and for all }x \in C_T.\]
For variables $u_{kl}$ in $I$ this is clear from the definition of $\sigma_T^{ij}$. Hence it holds in particular for all indices of the form $\{ik\}$ for $k \notin \{i,j\}$. In order to check this fact for the other indices, recall the definition of $I$ after Lemma 4.3.  Since $u_{j s_2} = u_{i s_1}\inv ( u_{s_1 s_2} + u_{i s_2} u_{j s_1})$, we find
\[\log |u_{j s_2}(\sigma_T^{ij}(x))| = \max\{-x_{i s_1} + x_{s_1 s_2}, \, - x_{i s_1} + x_{i s_2} + x_{j s_1} - x_{ij}\}.\]
Since $s_1$ and $s_2$ are in the same subtree (see figure \ref{tree}), we find that 
$x_{i s_1} + x_{j s_2} = x_{i s_2} + x_{j s_1} \geq x_{ij} + x_{s_1 s_2}$, 
which implies $\log |u_{j s_2}(\sigma_T^{ij}(x))| = x_{j s_2} - x_{ij}$. 
Inductively, we can show in this way our claim for all indices of the form $\{jk\}$. If we consider an index of the form $\{kl\}$ where $k,l \notin \{i,j\}$, we have $u_{kl} = u_{ik} u_{jl} - u_{il} u_{jk}$.  If $k$ is a leaf in the subtree $T_a$ and $l$ is a leaf in the subtree $T_b$ for $a <b$, we find $x_{kl} + x_{ij} = x_{il} + x_{jk} > x_{ik} + x_{jl}$. Hence the non-archimedean triangle inequality gives $\log|u_{kl}(\sigma_T^{ij}(x))| = \log|(u_{il} u_{jk})(\sigma_T^{ij}(x))| = x_{kl} -x_{ij}$. 
If $k$ and $l$ are leaves in the same subtree $T_a$, we may assume that $k \prec l$. If $k$ is the predecessor of $l$ in the total ordering $\preceq$ restricted to $T_a$, we are done, since then $u_{kl}$ is contained in $I$. If not, we let $m$ be the predecessor of $l$, so that $k \prec m \prec l$. 
The Pl\"ucker relations give $u_{im} u_{kl} = u_{ik} u_{ml} + u_{il} u_{km}$,
hence 
\[u_{kl} = u_{im}\inv (u_{ik}u_{ml} + u_{il} u_{km}).\]
Note that all variables on the right hand side except possibly $u_{km}$ are contained in $I$. Hence we can calculate $\log|u_{kl}(\sigma_T^{ij}(x))|$ by developing $u_{km}$ in a Laurent series in the variables in $I$. The resulting Laurent series does not allow any cancellation between the variables. Now we can apply the cherry property for the ordering $\preceq$, which says that 
$\{k,m\}$ or $\{m,l\}$ is a cherry for the quartet $\{i,k,m,l\}$. Hence we have 
\[x_{il} + x_{km} \leq x_{im} + x_{kl} = x_{ik} + x_{ml}\]
or 
\[x_{ik} + x_{ml} \leq x_{im} + x_{kl} = x_{il} + x_{km}.\]
In both cases, one can check directly that our claim holds.

In fact, the section $\sigma_T^{ij}$ has a more conceptual description. Let $\psi: \Gr_0(2,n) \rightarrow \mbox{Spec} K[u_{kl}^{\pm}: u_{kl} \in I] = \mathbb{G}_m^{2(n-2)}$ be the open embedding induced by (\ref{eq:birational}), and let $\Sigma$ be the standard skeleton of $(\mathbb{G}_m^{2(n-2)})^\an$ as in (\ref{eq:skeleton}). It is contained in the open analytic subvariety $\Gr_0(2,n)^\an$, since it consists entirely of norms and does therefore not meet a closed subvariety of strictly lower dimension. Then the map $\sigma_T^{ij}$ is  the composition of the projection from $\mathcal{C}_T$ to the coordinates $(x_{kl}- x_{ij})_{u_{kl} \in I} \in \mathbb{R}^{2(n-2)} = \Sigma$ followed by the inclusion of $\Sigma$ in $\Gr_0(2,n)^\an$.

Recall from section 2.3 that for every $x \in \cTG_0(2,n)$ the preimage $\trop_\varphi^{-1}(x)$ under the tropicalization map is the Berkovich spectrum $\mathcal{M}(A_x)$ of an affinoid algebra $A_x$. It follows from our construction that $\sigma_T^{ij}(x)$ is the unique Shilov boundary point of $A_x$. We can formulate this fact explicitely as follows. 

\begin{lem}\label{lem:maximal}[\cite{chw}, lemma 4.17]  For every $x \in \mathcal{C}_T$ the seminorm $\sigma_T^{ij}(x)$ constructed above has the following maximality property: For every $f$ in the coordinate ring $K[\Gr_0(2,n)]$ and every seminorm $p \in \trop_\varphi\inv (x) \subset \Gr_0(2,n)^\an$ we have 
\[|f(p)| \leq |f(\sigma_T^{ij}(x))|.\]
\end{lem}

As an immediate consequence we see that $\sigma_T^{ij}(x)$ does not depend on the choice of the partial ordering $\preceq$ or on the pair $ij$, and that $\sigma_{T}^{ij}(x) = \sigma_{T'}^{ij}(x)$ if $x$ is contained in the intersection $\mathcal{C}_T \cap \mathcal{C}_{T'}$. Hence we can patch these maps together and get a well-defined section $\sigma: \cTG_0(2,n) \rightarrow \Gr_0(2,n)^\an$ of the tropicalization map. It follows from the construction that it is continuous. 

This proves theorem \ref{thm:Grass} on the dense torus orbit $\Gr_0(2,n)$. In order to define a section of the tropicalization map also on the boundary strata of $\Gr(2,n)$, we follow the same strategy of constructing an index set $I$ such  that the associated Pl\"ucker variables generate the function field of the Grassmannian. This construction is more involved, see \cite{chw}, section 4 for details. The most important problem here is to find these local index sets in such a way that the section is also continuous when passing from one stratum of the Grassmannian for another. This continuity statement is shown in \cite{chw}, Theorem 4.19.

Using the index set $I$, we may also calculate the initial degenerations of all points in the tropical Grassmannian and deduce that their tropical multiplicity is one everywhere. We will explain a general argument showing the existence of a section in this case in paragraph \ref{subsection:copy}.  

In \cite{drpo}, Draisma and Postinghel give a different proof of Theorem \ref{thm:Grass}. They work with the affine cone over the Grassmannian $\Gr(2,n)$, define the section on a suitable subset and use torus actions and tropicalized torus actions to move it around. Also in this approach a maximality statement such as  Lemma \ref{lem:maximal} is used. 

\section{Skeleta of semistable pairs} 
In the next three sections, we give an overview of the results in \cite{grw}.
In order to compare polyhedral substructures of Berkovich spaces with tropicalizations, we start by generalizing Berkovich's notion of skeleta. Such skeleta are induced by the incidence complexes of the special fibers of suitable models. We extend this notion by adding a horizontal divisor on the model. For the rest of this paper, we fix an algebraically closed ground field $K$ which is complete with respect to a non-archimedean, non-trivial absolute value. 

\subsection{Integral affine structures} \label{section:integral} For curves, metrics play an important role in the comparison results between tropical and analytic varieties, as we have seen in section \ref{section:curves}. The right way to generalize this to higher dimensions is to consider integral affine structures. 
Let $M$ be a lattice in the finite-dimensional real vector space $M_{\mathbb{R}} =  M \otimes_{\mathbb{Z}} \mathbb{R}$, and let $N = \Hom(N,\mathbb{Z})$ its dual, which is a lattice in the dual space $N_{\mathbb{R}} = N \otimes_{\mathbb{Z}} \mathbb{R}$ of $M_\mathbb{R}$. 
We denote the associated pairing by $\langle \, , \, \rangle: M_\R \times N_\R \rightarrow \R$. Recall that $\Gamma = \log|K^\times|$. 

An \emph{integral $\Gamma$-affine polyhedron} in 
$N_\R$ is a subset of $N_\R$ of the form
\[ \Delta = \big\{ v\in N_\R \mid \langle u_i, v \rangle + \gamma_i \geq 0 
\mbox{ for all } i = 1,\ldots,r \big\} \]
for some $u_1,\ldots,u_r\in M$ and $\gamma_1,\ldots,\gamma_r\in\Gamma$.
Any face of an integral $\Gamma$-affine polyhedron $\Delta$ is again integral
$\Gamma$-affine. % The relative interior of $\Delta$ is denoted $\relint(\Delta)$.  
%A bounded polyhedron is called 
%a {\it polytope}.
 An \emph{integral $\Gamma$-affine polyhedral complex} in $N_\R$ is a polyhedral
complex whose faces are integral $\Gamma$-affine. 

An \emph{integral $\Gamma$-affine function} on $N_\R$ is a function from $N_\R$ to $\R$ which is  of the
form
\[ v \mapsto \langle u,v \rangle+ \gamma \]
for some $u\in M$ and $\gamma\in\Gamma$.  
More
generally, let $M'$ be a second finitely generated free abelian group and
let $N' = \Hom(M',\Z)$.  An \emph{integral $\Gamma$-affine map} from 
$N_\R$ to $N'_\R$ is a function of the form $F = \phi^* + v$, where
$\phi:M'\to M$ is a homomorphism, $\phi^*:N_\R\to N'_\R$ is the dual
homomorphism extended to $N_\R$, and $v\in N' \otimes_\mathbb{Z} \Gamma$.  If $N' = M' = \Z^m$
and $F = (F_1,\ldots,F_m) : N_\R\to\R^m$ is a function, then $F$ is
integral $\Gamma$-affine if and only if each coordinate $F_i:N_\R\to\R$ is
integral $\Gamma$-affine.  

An \emph{integral $\Gamma$-affine map} from an integral $\Gamma$-affine
polyhedron $\Delta\subset N_\R$ to $N'_\R$ is defined as the restriction to
$\Delta$ of an integral $\Gamma$-affine map $N_\R\to N'_\R$.  If
$\Delta'\subset N'_\R$ is an integral $\Gamma$-affine polyhedron then a
function $F:\Delta\to\Delta'$ is \emph{integral $\Gamma$-affine} if the
composition $\Delta\to\Delta'\hookrightarrow N'_\R$ is integral $\Gamma$-affine. We say that an integral $\Gamma$-affine map $F: \Delta \rightarrow N'_\R$ is \emph{unimodular} if $F$ is injective and if the inverse map $F(\Delta)\rightarrow \Delta$ is integral $\Gamma$-affine. Note that $F(\Delta)$ is an integral $\Gamma$-affine polyhedron in $N'_\R$.

\subsection{Semistable pairs}\label{section:semistable}
  Note that since our ground field $K$ is algebraically closed, the ring of integers $K^\circ$ is a valuation ring which is not discrete and not noetherian. 
We begin by decribing the building blocks of the polyhedral substructures of Berkovich spaces which lend themselves to comparison with tropicalizations.

Let $0 \leq r \leq d$ be natural numbers and consider the $K^\circ$-scheme
\[\mathscr{S} = \mbox{Spec}(K^\circ [x_0, \ldots, x_d] / (x_0 \ldots x_r - \pi))\] for some $\pi \in K^\circ$ satisfying $|\pi| < 1$. Then $\mathscr{S}$ is a flat scheme over $K^\circ$ with smooth generic fiber. Its special fibre contains $r+1$ irreducible components whose incidence complex is an $r$-dimensional simplex.
Now fix a natural number $s \geq 0$ such that $r+s \leq d$, and consider the principal Cartier divisor $H(s) = \mbox{div}(x_{r+1}) + \ldots + \mbox{div}(x_{r+s})$ on $\mathscr{S}$.  
%The pair $(\mathscr{S}, H(s))$ is called \emph{standard pair} in \cite{grw}, Example 3.10. 

As explained in section \ref{section:berkovich}, the $K^\circ$-scheme $\mathscr{S}$ gives rise to
two analytic spaces in the following way.
On the one hand, we have an associated admissible formal scheme 
\begin{equation}\label{eq:formalS}
\mathcal{S} = \mbox{Spf}\big(K^\circ \{ x_0,\ldots,x_d \}/ (x_0\ldots x_r-\pi )\big)
\end{equation}
which we get
by completing $\mathscr{S}$ with respect to a non-zero element in $K$ of absolute value $<1$. Its  analytic generic fiber is $\mathcal{M}(A_\mathcal{S})$ for the  affinoid algebra 
\[A_\mathcal{S}= K\{ x_0, \ldots, x_d \} / (x_0 \ldots x_r - \pi).\] 
It  is a subset of the analytification of the generic fibre $\mathscr{S}_K = \mbox{Spec}(K[x_0, \ldots, x_d]/ (x_0 \ldots x_r - \pi))$ of $\mathscr{S}$. 

Now we look at the tropicalization map on $\mathcal{M}(A_\mathcal{S})$ which only takes into account the first $r+s+1$ coordinates, i.e. the map
\begin{eqnarray*}
\Val:  \mathcal{M}(A_\mathcal{S}) \backslash \mathcal{H}& \rightarrow & \R^{r+s+1}_{\geq 0} \\
p & \mapsto & (-\log|x_0(p)|, \ldots, -\log|x_{r+s}(p)|),
\end{eqnarray*}
where $\mathcal{H}$ is the support of the Cartier divisor induced by $H(s)$. 
Its image is $\Delta(r,\pi) \times \mathbb{R}^s_{\geq 0}$, where $\Delta(r,\pi) = \{(v_0, \ldots, v_r) \in \mathbb{R}_{\geq 0}^{r+1}: v_0 + \ldots + v_r = -\log|\pi|\}$ is a simplex in $\R^{r+1}$.

We define a continuous section $\Delta(r,\pi) \times \mathbb{R}^s_{\geq 0} \rightarrow \mathcal{M}(A_\mathcal{S})$ as follows. 
Note that the projection  $(x_0, x_1, \ldots, x_d) \mapsto (x_1, \ldots, x_d)$ induces an isomorphism from
$\mathcal{M}(A_\mathcal{S})$ to the affinoid subdomain $B = \{p: - \log|x_1(p)| - \ldots - \log|x_r(p)| \leq -\log|\pi|\}$ of $\mathcal{M}(K\{ x_1, \ldots, x_d \})$. 

Similarly, the projection  $(v_0, v_1, \ldots, v_{r+s}) \mapsto (v_1, \ldots, v_{r+s})$ induces a homeomorphism \[\Delta(r,\pi) \times \R^s_{\geq 0} \rightarrow \Sigma = \{(v_1, \ldots, v_{r+s}) \in \R_{\geq 0}^{r+s}: v_1 + \ldots + v_r \leq  -\log|\pi|\}\]

For every $v = (v_1, \ldots, v_{r+s} ) \in \Sigma$ there is a bounded multiplicative norm $|| \,\, ||_v$ on $K\{ x_1, \ldots, x_d \}$ which is defined as follows:
\[ || \sum_{I = (i_1, \ldots, i_d)} a_I x^I||_v = \max_I \{|a_I| \exp(-i_1 v_1 - \ldots - i_{r+s} v_{r+s})\}.\]
It satisfies $- \log||x_1||_v - \ldots - \log||x_r||_v = v_1 + \ldots + v_r \leq  - \log |\pi|$. Therefore the point $|| \,\,||_v$ is contained in the affinoid domain $B$.  Hence there is a uniquely determined continuous map 
\[\sigma: \Delta(r,\pi) \times \R_{\geq 0}^{s} \rightarrow  \mathcal{M}(A_\mathcal{S})\]
making the diagram
\[
\begin{xy}

\xymatrix{
\Delta(r,\pi)  
\times \R_{\geq 0 }^s  
\ar[r]^\sigma \ar[d] & \mathcal{M}(A_\mathcal{S})  \ar[d] \\
\Sigma  \ar[r]^{v \mapsto || \, \, ||_v} & B 
 }
\end{xy}
\]
commutative. 
The map $\sigma$ is by construction a section of the map $\mbox{Val}$. 
We define $S(\mathcal{S}, H(s))  \subset \mathcal{M}(A_\mathcal{S})  \subset (\mathscr{S}_K)^\an$ as the image of $\sigma$ and call it the skeleton of the pair $(\mathcal{S}, H(s))$. 

Now we consider schemes which \'etale locally look like some $\mathscr{S}$.

\begin{defi} 
\label{strictly semistable pairs}
  A {\it strictly semistable pair} $(\mathscr{X}, H)$ consists of an irreducible proper flat
  scheme $\mathscr{X}$ over the valuation ring $K^\circ$ and a sum
  $H = H_1 + \cdots + H_S$ of  effective Cartier divisors
  $H_i$ on $\mathscr{X}$ such that $\mathscr{X}$ is covered by open  
  subsets $\mathscr{U}$ which admit an \'etale morphism 
  \begin{equation}
    \label{eq:local.sss.pair}
    \psi ~:~ {\mathscr{U}} \longrightarrow {\mathscr{S}} =
    \mbox{Spec}( K^\circ[ x_0, \dots, x_d ] / (x_0 \cdots x_r - \pi) ) 
  \end{equation}
  for some $r \leq d$ and $\pi \in K^\times$ with $|\pi|<1$.  We assume that
  each  $H_i$ has irreducible support and that the restriction of $H_i$ to $\mathscr{U}$ is either trivial or defined by 
  $\psi^*(x_{j})$ for some $j \in \{r+1, \dots, d\}$.
\end{defi}

This is a generalization of de Jong's notion of a strictly semistable pair over a discrete valuation ring \cite{jo}, if we add the divisor of the special fiber to the horizontal divisor $H$. It is sometimes convenient to include only the horizontal divisor $H$ as part of the data, since it is a Cartier divisor, whereas the special fiber of $\mathscr{X}$ may not be one. 
For a more detailed discussion of this issue see \cite{grw}, Proposition 4.17.
If $d=1$, i.e. in the case of curves, semistable models and therefore nice skeletons are always available after an extension of the ground field. In higher dimensions, we have to allow alterations in order to find semistable models, see \cite{jo}, Theorem 6.5.

\subsection{Skeleta} \label{section:skeleta} Let $(\mathscr{X},H)$ be a semistable pair in the sense of definition \ref{strictly semistable pairs}. The generic fiber $X$ of $\mathscr{X}$  is smooth of dimension $d$. Hence $d$ is constant in every chart $\mathscr{U}$, whereas the numbers $r$ and $s$  may vary with $\mathscr{U}$. We denote the special fiber of $\mathscr{X}$ by $\mathscr{X}_s= \mathscr{X} \otimes_{K^\circ} \tilde{K}$. 

Denote by $V_1, \ldots, V_R$ the irreducible components of the special fiber $\mathscr{X}_s$ of $\mathscr{X}$. The Cartier divisors $H_i$ which are part of our data give rise to horizontal closed subschemes $\mathscr{H}_i$ of $\mathscr{X}$, which are locally cut out by a defining equation of $H_i$. Putting $D_i = V_i$ for $i = 1, \ldots, R$ and $D_{i + R} = \mathscr{H}_i$, we get a Weil divisor $D = \sum_{i = 1}^{R+S} D_i$ on $\mathscr{X}$. This gives rise to a stratification of $\mathscr{X}$, where a stratum is defined as an irreducible component of a set of the form $\bigcap_{i \in I} D_i \backslash \bigcup_{i \notin I} D_i$ for some $I \subset \{1, \ldots, R+S\}$. We call any stratum contained in the special fiber $\mathscr{X}_s$ a vertical stratum and denote the set of all vertical strata by $\mbox{str}(\mathscr{X}_s, H)$. 

Now we want to glue skeleta of local charts together in such a way that the faces of the resulting polyhedral complex are in bijective correspondence with the  vertical strata  in $\mbox{str}(\mathscr{X}_s, H)$. In order to achieve this, we may have to pass to a smaller covering in the category of admissible formal schemes.  

Let $(\mathscr{X},H)$ be a strictly semistable pair with a covering as in Definition \ref{strictly semistable pairs}. We consider the induced formal open covering of the associated admissible formal scheme $\mathcal{X}$ which is defined by completion. Hence $\mathcal{X}$ is covered by formal open subsets $\mathcal{U}$ which admit an \'etale morphism $\psi: \mathcal{U} \rightarrow \mathcal{S}$ to a formal scheme $\mathcal{S}$ as in (\ref{eq:formalS}).

It is shown in \cite{grw}, Proposition 4.1 that, after passing to a refinement,  the formal \'etale covering 
\[\psi: \mathcal{U} \rightarrow \mathcal{S}\]
has the property that $\psi\inv\{x_0 = \ldots = x_{r+s} = 0\}$ is a vertical stratum $S$ in the special fiber $\mathscr{X}_s$ such that for every  vertical stratum $T$ the following condition holds: The closure $\overline{T}$ of $T$ in $\mathscr{X}_s$ meets $\mathcal{U}_s$ if and only if $S \subset \overline{T}$. 

We define the skeleton of $(\mathcal{U}, H|_{\mathcal{U}})$ as the preimage of the skeleton of the standard pair: $S(\mathcal{U}, H|_{\mathcal{U}}) =  \psi^{-1} (S(\mathcal{S}, H(s))$. This is a subset of the analytic generic fiber $\mathcal{U}_\eta$ of $\mathcal{U}$ which does not meet the horizontal divisor $H$. 

It follows from results of Berkovich \cite{ber99} that the \'etale map $\psi$ actually induces a homeomorphism between
$S(\mathcal{U}, H|_{\mathcal{U}})$ and $S(\mathcal{S}, H(s)) \simeq \Delta(\pi,r) \times \R_{\geq 0}^{s}$. 
Hence the map 
\[\Val \circ \psi: S(\mathcal{U}, H|_{\mathcal{U}}) \rightarrow \Delta(\pi,r) \times \R_{\geq 0}^s\]
is a homeomorphism.

It is shown in \cite{grw}, 4.5 that the skeleton $S(\mathcal{U}, H|_{\mathcal{U}})$ only depends on the minimal stratum $S$ contained in the special fiber of $\mathcal{U}$. Therefore we denote it by 
 $\Delta_S$ and call it the \emph{canonical polyhedron of $S$}. The dimensions of $S$ and of the  canonical polyhedron $\Delta_S$ are defined in an obvious way and add up to $d$, see \cite{grw}, Proposition 4.10.
 
 As we have seen, $\Delta_S$ is homeomorphic to $\Delta(\pi,r) \times \R_{\geq 0}^{s}$ for suitable data $r,s$ and $\pi$ as above. We call $\Delta(\pi,r)$ the \emph{finite part} and $\R_{\geq 0}^s$ the \emph{infinite part} of $\Delta_S$.
 \begin{figure}[hbf]
  \centering
  \begin{minipage}{0.3\textwidth}
    \includegraphics[width=0.3\textwidth, height=1.5cm]{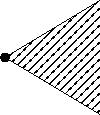}
    \end{minipage}
    \begin{minipage}{0.3\textwidth}
    \includegraphics[width=0.4\textwidth, height=1cm]{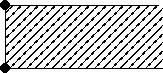}
    \end{minipage}
    \begin{minipage}{0.3\textwidth}
    \includegraphics[width=0.3\textwidth, height=1.2cm]{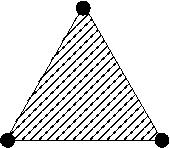}
  \end{minipage}
\caption{The canonical polyhedron $\Delta_S$ in the case $r=0$ and $s=2$, in the case $r=s=1$ and in the case $r=2$ and $s=0$ (from left to right).}
\label{finitepart}
\end{figure}
 
  As explained above, $S \mapsto \Delta_S$ is a bijective correspondence between faces $\Delta_S$ of the skeleton $S(\mathscr{X}, H)$ and vertical strata $S$ induced by the divisor $D$ in the special fiber $\mathscr{X}_s$. Recall that $D$ is given by the horizontal divisor $H$ plus all irreducible components of $\mathscr{X}$ (which we also call vertical divisors). Note that a vertical stratum $T$ satisfies $S \subset \overline{T}$ if and only if $\Delta_T$ is a closed face of $\Delta_S$.

Now we can glue $\Delta_S$ and $\Delta_{S'}$ along the union of all canonical polyhedra associated to vertical strata $R$ such that $\overline{R} \supset S \cup S'$. 
In this way we define the skeleton of $(\mathscr{X}, H)$ as the union of all $\Delta_S$ for strata $S$ in the special fiber $\mathscr{X}_s$:
\[S(\mathscr{X}, H) = \bigcup_{S} \Delta_S.\]
We obtain a piecewise linear space $S(\mathscr{X}, H)$  whose charts are integral $\Gamma$-affine polyhedra.

The skeleton $S(\mathscr{X},H)$ is a closed subset of the analytic space $X^\an \backslash H_K$, where we  write $H_K$ for the generic fiber of the support of $H$ on $\mathscr{X}$. Using methods from \cite{ber99}, one can prove that it is in fact a strong deformation retract of $X^\an \backslash H_K$. In \cite{grw}, Theorem 4.13 it is shown that the retraction map can be extended to a retraction map from $X^\an$ to a suitable compactification of the skeleton. 

If the horizontal divisor $H = 0$, then the skeleton $S(\mathscr{X},0)$ is equal to Berkovich's skeleton of a semistable scheme, see \cite{ber99}.

\section{Functions on the skeleton}
Throughout this section we  fix a strictly semistable pair $(\mathscr{X}, H)$. We use the notation from the previous sections. 

We want to show that for every non-zero rational function $f$ on $X$ such that the support of the divisor of $f$ is contained in $H_K$, the function $- \log |f|$ factors through a piecewise integral $\Gamma$-affine map on the skeleton. Moreover, we will show a slope formula for this map. 
 
\begin{thm}[\cite{grw}, Proposition 5.2] \label{function} Let $f$ be a non-zero rational function on $X$ such that the support of $\mbox{div}(f)$ is containd in $H_K$. We put $U = X \backslash H_K$ and consider the function 
\[F = - \log|f|: U^{\an} \rightarrow \R.\]
Then $F$ factors  through the retraction map $\tau: U^{\an} \rightarrow S(\mathscr{X}, H)$ to the skeleton:

\[ \xymatrix{
U^\an \ar[rr]^{F} \ar[rd]_{\tau}& &\mathbb R \\
&S(\mathscr{X},H)  \ar[ru]_{F|_{S(\mathscr{X},H)}}&
} \]

Moreover, the restriction of $F$ to $S(\mathscr{X},H)$ is an integral $\Gamma$-affine function on each canonical polyhedron.
\end{thm}

The theorem is proved by considering the formal building blocks with distinguished strata and using a result of Gubler \cite{gu07}, Proposition 2.11.

Recall that we denote the dimension of $X$ by $d$. 
The slope formula for $F$ is basically a balancing condition around each $(d-1)$-dimensional canonical polyhedron of the skeleton which involves slopes in the direction of all adjacent $d$-dimensional polyhedra. 

Let $\Delta_S$ be a $d$-dimensional canonical polyhedron of the skeleton $S(\mathscr{X},H)$ containing the $(d-1)$-dimensional canonical polyhedron $\Delta_T$. Then   the stratum $S$ in the special fiber $\mathscr{X}_s$ is contained in the closure $\overline{T}$ (which is a curve), and  it is obtained as a component of the intersection of $\overline{T}$ with one additional irreducible component of the divisor $D$.

 This component is either vertical, i.e. a component of the special fiber,  or horizontal, i.e. given by some component of $H$. If it is vertical, then the finite part of $\Delta_S \simeq \Delta(r,\pi) \times \R^{s}_{\geq 0}$ (i.e. the simplex $\Delta(r,\pi)$) is strictly larger than the finite part of $\Delta_T$. In this case we say that $\Delta_S$ \emph{extends $\Delta_T$ in a bounded direction}.  If the component is horizontal, then the infinite part of $\Delta_S$ (i.e. the product $\R^s_{\geq 0}$) is strictly larger than the infinite part of $\Delta_T$. In this case we say that $\Delta_S$ \emph{extends $\Delta_T$ in an unbounded direction}.
  
We define the \emph{bounded  degree} $\degree_b(\Delta_T)$ as the number of canonical polyhedral $\Delta_S$ extending $\Delta_T$ in a bounded direction. Similarly, we define the \emph{unbounded degree} $\degree_u(\Delta_T)$ as the number of canonical polyhedra $\Delta_S$ extending $\Delta_T$ in an unbounded direction. 

In the case $d=1$, i.e. if $X$ is a curve,  all $\Delta_T$ are vertices. A one-dimensional canonical polyhedron $\Delta_S$ extending $\Delta_T$ in a bounded direction is simply an edge of length $v(\pi) > 0$. In this case, the stratum $S$ is a component of the intersection of two irreducible components in the special fiber $\mathscr{X}_s$. A one-dimensional canonical polyhedron $\Delta_S$ extending $\Delta_T$ in an unbounded direction is a ray of the form $\R_{\geq 0}$. In this case, the  stratum $S$ is a component of the intersection of the irreducible component of the special fiber given by $T$ with a horizontal component of the divisor. Hence the degree $\degree_b(\Delta_T)$ is the number of bounded edges in $\Delta_T$, and $\degree_u(\Delta_T)$ is the number of unbounded rays starting in the vertex $\Delta_T$. 

Now we want to define multiplicities via intersection theory. Since $\mathscr{X}$ is not noetherian, the standard intersection theory tools from algebraic geometry are not available. However, on admissible formal schemes one can use analytic geometry to associate Weil divisors to Cartier divisors, and there is a refined intersection product with Cartier divisors, see \cite{gu98}, \cite{gu03} and Appendix A in \cite{grw}. 

\begin{defi} For every vertex $u \in \Delta_T$ we denote by $V_u$ the associated irreducible component of the special fiber $\mathscr{X}_s$. We define an integer $\alpha(u, \Delta_T)$ as follows:

If $\Delta_T$ has zero-dimensional finite part $\{u\}$, we simply put $\alpha(u, \Delta_T) = \degree_b(\Delta_T)$.

 If not, then $T$ lies in at least two irreducible components of the special fiber $\mathscr{X}_s$, which means that the finite part of $T$ is a simplex $\Delta(r,\pi)$ of dimension at least one. If $\mathcal{U}$ is a suitable formal chart, it is shown in \cite{grw}, Proposition 4.17 that there exists a unique effective Cartier divisor $C_u$ on $\mathcal{U}$ such that its Weil divisor is equal to $v(\pi) (V_u \cap \mathcal{U}_s)$. We define $-\alpha(u, \Delta_T)$ as the intersection number $ (C_u . \overline{T})$, which is equal to the degree of the pullback of the line bundle associated to $C_u$ from $\mathcal{U}$ to $\overline{T}$. 

\end{defi}

Note that Cartwright \cite{ca} introduced the intersection numbers $\alpha(u, \Delta_T)$ in the (noetherian) situation where $K^\circ$ is discrete valuation ring. He uses them to endow the compact skeleton $S(\mathscr{X})$ with the structure of a tropical complex, see \cite{ca}, Definition 1.1. This notion also involves a local Hodge condition which plays no role in our slope formula. The paper \cite{ca} develops a theory of divisors on tropical complexes and investigates their relation to algebraic divisors.

Here we also need multiplicities for rays in $\Delta_T$, which are defined as one-dimensional faces of the unbounded part of $T$.

\begin{defi}Let $H_r$ be the horizontal component corresponding to the ray $r$ in $\Delta_T$. We put
\[\alpha(r,\Delta_T) = - (H_r.\overline{T}),\]
where we take the intersection product of $\overline{T}$ with the Cartier divisor $H_r$ on $\mathscr{X}$.
\end{defi}

Let $F$ be a function on the skeleton $S(\mathscr{X},H)$, which is  integral $\Gamma$-affine on each canonical polyhedron.
Now we are ready to define outgoing slopes of $F$  on a $(d-1)$-dimensional canonical polyhedron $\Delta_T$ along a $d$-dimensional polyhedron $\Delta_S$. 

\begin{defi}
Let $\Delta_T$ be a $(d-1)$-dimensional canonical polyhedron and let $\Delta_S$ be a $d$-dimensional canonical polyhedron of $S(\mathscr{X},H)$ containing $\Delta_T$. We denote by $\Delta(r,\pi)$  the finite part of $\Delta_S$. Let $F: \Delta_S \rightarrow \mathbb{R}$ be an integral $\Gamma$-affine function.

i) If $\Delta_S$ extends $\Delta_T$ in a bounded direction, then there exists a unique vertex $w$ of $\Delta_S$ not contained in $\Delta_T$. We put 
\[\slope(F; \Delta_T, \Delta_S) = \frac{1}{v(\pi)} \left(F(w) - \frac{1}{\degree_b(\Delta_T)} \sum_{u \in \Delta_T} \alpha(u, \Delta_T) F(u) \right),\]
where we sum over all vertices $u$ in $\Delta_T$. 

ii) If $\Delta_S$ extends $\Delta_T$ in an unbounded direction, then there exists a unique ray $s$ in $\Delta_S$ not contained in $\Delta_T$, and we put
\[\slope(F; \Delta_T, \Delta_S) = d_s F - \frac{1}{\degree_u(\Delta_T)} \sum_{r \in \Delta_T} \alpha(r, \Delta_T) d_r F ,\]

where we sum over all rays in $\Delta_T$. For any ray $r$ we denote by $d_r F$ the derivative of $F$ along the primitive vector in the direction of $r$. 
\end{defi}

If $X$ is a curve, $\Delta_T= u$ is a vertex and $\Delta_S$ is an edge with vertices $u$ and $w$, then  $\slope(F; \Delta_T, \Delta_S) = \frac{1}{v(\pi)} (F(w) - F(u))$.
If $\Delta_S$ is a ray $s$ starting in $u$, then we simply have $\slope(F; \Delta_T, \Delta_S) = d_s F $.
In higher dimensions, the definition is more involved, since the naive slope 
$\frac{1}{v(\pi)} (F(w) - F(u))$ depends on the choice of a vertex $u$ in $\Delta_T$. Therefore we define a replacement for a weighted midpoint in $\Delta_T$ as 
$\frac{1}{\degree_b(\Delta_T)} \sum_{u \in \Delta_T} \alpha(u, \Delta_T) u$. Note that this point does not necessarily lie in $\Delta_T$, since $\alpha(u, \Delta_T)$ may be negative. 

We can now formulate the slope formula for skeleta.

\begin{thm}[\cite{grw}, Theorem 6.9]\label{slope_formula} Let $f \in K(X)^\times$ be a non-zero rational function such that the support of $\mbox{div}(f)$ is contained in $H_K$. Let $F: S(\mathscr{X},H) \rightarrow \R$ be the restriction of the function $ -\log|f|$ to the skeleton. Then $F$ is continuous and integral $\Gamma$-affine on each canonical polyhedron of $S(\mathscr{X},H)$, and for all $(d-1)$-dimensional canonical polyhedra we have
\[\sum_{\Delta_S \succ \Delta_T} \slope(F; \Delta_T, \Delta_S) = 0,\]
where the sum runs over all $d$-dimensional canonical polyhedra $\Delta_S$ containing $\Delta_T$. 
\end{thm}

If $X$ is a curve, the slope formula basically says that the sum of all outgoing slopes along edges or rays in a fixed vertex  is zero. In this case, the slope formula is shown in \cite[Theorem 5.15]{bpr2}. It is a reformulation of the non-archimedean Poinca\'re-Lelong formula proven in Thuillier's thesis \cite{thu}, Proposition 3.3.15. The Poincar\'e-Lelong formula is an equation of currents in the form $d d^c \log|f| = \delta_{\mathrm{div}(f)}$, where $d d^c$ is a certain distribution-valued operator. This version of the slope formula was generalized to higher dimensions in the ground-breaking paper
\cite{cldu}, where a theory of differential forms and currents on Berkovich spaces is developed. The approach of \cite{cldu} uses tropical charts and does not rely on models or skeleta. In higher dimensions we see no direct relation to our slope formula.

\section{Faithful tropicalizations}\label{section:faithful}

We will now investigate the relation between skeleta, which are polyhedral substructures of analytic varieties, and tropicalizations, which are polyhedral images of algebraic or analytic varieties. 

\subsection{Finding a faithful tropicalization for a skeleton}
We start with a strictly semistable pair $(\mathscr{X}, H)$ with skeleton $S(\mathscr{X},H)$ and generic fiber $X$. We consider rational maps $f:X \dashrightarrow \mathbb{G}_{m,K}^n$ from $X$ to a split torus. If $U \subset X$ is a Zariski open subvariety where $f$ is defined, then $f|_U: U \rightarrow \mathbb{G}_{m,K}^n$ induces a tropicalization $\Trop_{f}(U)$ of $U$, which is defined as the image of the map
$\trop \circ f^{\an}: U^{\an} \rightarrow \R^n$ as in section \ref{section:tropicalization}.

Note that the skeleton is contained in the analytification of every Zariski open subset of $X$, since it only contains norms on the function field of $X$. 

\begin{defi} A rational map $f:X \dashrightarrow \mathbb{G}_{m,K}^n$ from $X$ to a split torus is called a faithful tropicalization of the skeleton $S(\mathscr{X},H)$ if the following conditions  hold:

i) The map $\trop \circ f^{\an}$ is injective on $S(\mathscr{X},H)$.  

ii) Each canonical polyhedron $\Delta_S$ of $S(\mathscr{X},H)$ can be covered by finitely many integral $\Gamma$-affine polyhedra such that  the restriction of $\trop \circ f^{\an}$ to each of those polyhedra is a unimodular integral $\Gamma$-affine map. 

\end{defi}

For the definition of unimodular integral affine maps see section \ref{section:integral}.

It is easy to see that a rational map which is unimodular on the skeleton stays unimodular if we enlarge it with more rational functions on $X$, see \cite{grw}, Lemma 9.3.

If we look at a building block $\mathscr{S} = \mbox{Spec} K^\circ [x_0, \ldots, x_d] / (x_1 \ldots x_r - \pi)$ as in Definition \ref{strictly semistable pairs} with the local tropicalization $\Delta(r,\pi) \times \R^s_{\geq 0}$, we find that the coordinate functions $x_0, \ldots, x_{r+s}$ induce a faithful tropicalization. Collecting the corresponding rational functions on $X$ for all canonical polyhedra of the skeleton, we get a rational map on $X$ which is locally unimodular on the skeleton. It is shown in the proof of Theorem 9.5 of \cite{grw} how to enlarge this collection of rational function in order to ensure injectivity of the tropicalization map on the skeleton. In this way one can show the following result.

\begin{thm}[\cite{grw}, Theorem 9.5] Let $(\mathscr{X},H)$ be a strictly semistable pair. Then there exists a collection of non-zero rational functions $f_1, \ldots, f_n$ on $X$ such that the resulting rational map $f=(f_1, \ldots, f_n): X \dashrightarrow \mathbb{G}^{n}_{m,K}$ is a faithful tropicalization of the skeleton $S(\mathscr{X},H)$. 
\end{thm}

\subsection{Finding a copy of the tropicalization inside the analytic space} \label{subsection:copy}
Let us now start with a given tropicalization of a very affine $K$-variety $U$, i.e.
with a closed immersion $\varphi: U \hookrightarrow \mathbb{G}_{m,K}^n$ of a variety $U$ in a split torus. As in section \ref{section:tropicalization} we  consider the tropicalization map
\[\trop_\varphi = \trop \circ \varphi^{\an}: U^{\an} \rightarrow \mathbb{R}^n.\]

Recall  that for every point $\omega \in \Trop_\varphi(U)$ the preimage $\trop_\varphi^{-1}(\omega) \subset U^{\an}$ of the tropicalization map is the Berkovich spectrum of an affinoid algebra $A_\omega$. 

\begin{lem}[\cite{grw}, Lemma 10.2]  \label{shilov} If the tropical multiplicity of $\omega$ is equal to one, then $A_\omega$ contains a unique Shilov boundary point.
\end{lem}

This lemma follows basically from \cite{bpr}, Remark after Proposition 4.17, which shows how the preimage of tropicalization is related to initial degenerations.

Now we can show that on the locus of tropical multiplicity one there exists a natural section of the tropicalization map. 

\begin{thm} [\cite{grw}, Theorem 10.7]\label{section} Let $Z \subset \Trop_\varphi(U)$ be a subset such that the tropical multiplicity of every point in $Z$ is equal to one. By Lemma \ref{shilov} this implies that for every $\omega \in Z$ the affinoid space $\trop_\varphi^{-1}({\omega})$ has a unique Shilov boundary point which we denote by $s(\omega)$. 

The map $s: Z \rightarrow U^{\an}$, given by $ \omega \rightarrow s(\omega)$, is continuous and a partial section of the tropicalization map, i.e. on $Z$ we have $\trop_\varphi \circ s = \mbox{id}_Z$. 
Hence the image $s(Z)$ is a subset of $U^{\an}$ which is homeomorphic to $Z$.

Moreover, if $Z$ is contained in the closure of its interior in $\Trop_\varphi(U)$, then $s$ is the unique continuous section of the tropicalization map on $Z$.
 \end{thm}

We give a sketch of the proof. 
It is enough to show that the section $s$ is continuous and uniquely determined under our additional assumption. Since everything behaves nicely under base change we may assume that the valuation map $K^\times \rightarrow \R_{>0}$ is surjective.  Let us first consider the case that $U = \mathbb{G}_{m,K}^n$ and $\varphi$ the identity map. Then the section $s$ is the identification of the tropicalization $\mathbb{R}^n$ (which has multiplicity one everywhere) with the skeleton of  $\mathbb{G}_{m,K}^n$ as defined in (\ref{eq:skeleton}). It is clear from the explicit description of $s$ in (\ref{eq:skeleton}) that it is continuous in this case.

 Moreover, if $s'$ is a different section of the tropicalization map in the case $U = \mathbb{G}_{m,K}^n$ which satisfies $s(\omega) \neq s'(\omega)$, we find a Laurent polynomial $f$ on which those two seminorms differ. Since $s(\omega)$ is the unique Shilov boundary point in the fiber of tropicalization, $s'(\omega)$ applied to $f$ is strictly smaller than $s(\omega)$ applied to $f$. Since $s(\omega) $ is equal to $s'(\omega)$ on all monomials, it follows that the initial degeneration of $f$ at $\omega$ cannot be a monomial. 
Therefore $\omega$ is contained in the tropical hypersurface $\Trop(f)$. A continuity argument shows that the same argument works in a small neighbourhood of $\omega$. This is a contradiction  since $\Trop(f)$ has codimension one. Therefore $s$ is indeed uniquely determined if $\varphi$ is the identity map. 

For general $\varphi: U \hookrightarrow \mathbb{G}_{m,K}^n$, where $U$ has dimension $d \leq n$, one shows that that there exists a  linear map $\Z^n \rightarrow \mathbb{Z}^d$ with the following property:
Let $\alpha: \mathbb{G}^n_{m,K} \rightarrow \mathbb{G}^d_{m,K}$ be the corresponding homomorphism of tori and consider  $\psi = \alpha \circ \varphi: U \rightarrow \mathbb{G}_{m,K}^d$. Let $S(\mathbb{G}^d_{m,K} )$ be the skeleton of the torus as in (\ref{eq:skeleton}).
Then for all $\omega \in Z$  we have $\{s(\omega)\} = \trop_\varphi\inv (\omega) \cap \psi\inv(S(\mathbb{G}^d_{m,K}))$. This implies that $s(Z) = \trop_\varphi\inv(Z) \cap \psi\inv(S(\mathbb{G}^d_{m,K}))$ is closed, from which we can deduce continuity of $s$. Uniqueness follows from uniqueness in the torus case by composing the sections with $\psi$. 

Note that the preceeding theorem does not make any assumption on the existence of specific models. If we assume that $(\mathscr{X}, H)$ is a semistable pair and $U = X \backslash H_K$, then the image of the section $s$ defined in Theorem \ref{section} is contained in the skeleton $S(\mathscr{X}, H)$. This is shown in \cite{grw}, Proposition 10.9.

The preceeding theorem treats the case of tropicalizations in tori. Let $\varphi: X \hookrightarrow Y$ be a closed embedding of $X$ in a toric variety $Y$ associated to the fan $\Delta$. Then we may consider the associated tropicalization $\Trop_\varphi(X)$ of $X$, i.e. the image of $\trop \circ \varphi^\an : X^\an \rightarrow Y^\an \rightarrow N_\mathbb{R}^\Delta$, where $N_{\mathbb{R}}^\Delta$ is the associated partial compactification of $N_\R$, see section \ref{section:tropicalization}. We can apply Theorem \ref{section} to all torus orbits. In this way, we get a section of the tropicalization map on the locus $Z \subset \Trop_\varphi(X)$ of tropical multiplicity one which is continuous on the intersection with each toric stratum. 
It is a natural question under which conditions this section is continuous on the whole of $Z$. This might shed new light on Theorem \ref{thm:Grass} for the tropical Grassmannian.

\begin{center}
Institut f\"ur Mathematik\\
Goethe-Universit\"at Frankfurt\\
Robert-Mayer-Strasse 8\\
D- 60325 Frankfurt\\
email: werner@math.uni-frankfurt.de
\end{center}

\end{document}